\documentclass[11pt,a4paper]{article}
\usepackage{ifthen,latexsym,amssymb,amsmath,bbm,fixmath}
\usepackage[nobysame,initials]{amsrefs}
\usepackage[shortlabels]{enumitem}
\usepackage{todonotes}
\newcommand{\eps}{\varepsilon}
\newcommand{\C}[1]{{\protect\mathcal{#1}}}

\newcommand{\I}[1]{{\mathbbm #1}}

\newcommand{\me}{{\mathrm e}}
\renewcommand{\mid}{:}
\renewcommand{\ldots}{\hspace{0.9pt}.\hspace{0.3pt}.\hspace{0.3pt}.\hspace{1.5pt}}
\renewcommand{\ge}{\geqslant}
\renewcommand{\le}{\leqslant}
\renewcommand{\geq}{\geqslant}
\renewcommand{\leq}{\leqslant}

\newcommand{\actson}{{\curvearrowright}}
\newcommand{\dd}{\,\mathrm{d}}

\newif\ifnotesw\noteswtrue

\newcommand{\hide}[1]{}

\newcommand{\beq}[1]{\begin{equation}\label{#1}}
	\newcommand{\eeq}{\end{equation}}

\newtheorem{theorem}{Theorem}
\newtheorem{lemma}[theorem]{Lemma}

\newcommand{\bpf}[1][Proof.]{\smallskip\noindent{\it #1} }
\newcommand{\qed}{\nolinebreak\mbox{\hspace{5 true pt}%
		\rule[-0.85 true pt]{3.9 true pt}{8.1 true pt}}}
\newcommand{\epf}{\qed \medskip}

\newtheorem{claim}{Claim}[theorem]

\newcommand{\SO}{\mathrm{SO}}
\newcommand{\Ra}[1]{\cite[#1]{Ratcliffe19fhm}}
\newcommand{\BR}[1]{\cite[#1]{BowenRadin03}}

\newcommand{\HMeasure}{\mu_m}
\newcommand{\vol}{\HMeasure}
\newcommand{\CO}{\C O}
\newcommand{\CG}{\C G}
\newcommand{\CL}{\C L}
\newcommand{\dist}{d}
\newcommand{\X}{X}
\newcommand{\Y}{Y}
\newcommand{\Orbit}{O}
\newcommand{\Ball}{B}
\newcommand{\D}{D}

\title{New lower bound on ball packing density\\ in high-dimensional hyperbolic spaces}

\author{
	Irene Gil Fern\'andez
	\thanks{Mathematics Institute and DIMAP, University of Warwick, UK. Email: {\tt irene.gil-fernandez@warwick.ac.uk}. IGF was supported by the Warwick Mathematics Institute Centre for Doctoral Training, and gratefully acknowledges funding from the University of Warwick and the UK Engineering and Physical Sciences Research Council (Grant number: EP/TS1794X/1)}
	\and
	Jaehoon Kim
	\thanks{Department of Mathematical Sciences, KAIST, South Korea. Email: {\tt jaehoon.kim@kaist.ac.kr}. J.K. was supported by the National Research Foundation of Korea (NRF) grant funded by the
	Korea government(MSIT) No. RS-2023-00210430}
	\and
	Hong Liu
	\thanks{Extremal Combinatorics and Probability Group (ECOPRO), Institute for Basic Science (IBS), Daejeon, South Korea. Email: {\tt hongliu@ibs.re.kr}. H.L. was supported by the Institute for Basic Science (IBS-R029-C4).}
	\and
	Oleg Pikhurko
	\thanks{Mathematics Institute and DIMAP, University of Warwick, UK.  Email: {\tt O.Pikhurko@warwick.ac.uk}. O.P. was supported by ERC Advanced Grant 101020255.}
}

\begin{document}

\maketitle	

\begin{abstract}
We present a new lower bound on the Bowen--Radin maximal density of radius-$R$ ball packings in the $m$-dimensional hyperbolic space, improving on the basic covering bound by factor $\Omega(m(R+\ln m))$ as $m$ tends to infinity. This is done by applying the recent theorem of Campos, Jenssen, Michelen and Sahasrabudhe 
%\cite{CamposJenssenMichelenSahasrabudhe} 
%who achieved an improvement of factor $\Omega(m\ln m)$ for the $m$-dimensional Eucledian space.
on independent sets in graphs with sparse neighbourhoods.
\end{abstract}

\section{Introduction}

Let $R>0$ be a positive real and let $(V,\dist)$ be a metric space endowed with a Borel non-zero measure~$\mu$. 
An \emph{$R$-packing} in $V$
%in a metric space $(V,\dist)$ 
is a subset $\X\subseteq V$ such that any two distinct points of $\X$ are at distance at least~$2R$. If the measure $\mu$ is \emph{finite} (that is, $\mu(V)<\infty$), then we define the \emph{density} of an $R$-packing $\X$ by 
 $$
 \D_R(\X):=
 \frac{\mu\left(\Ball_R(\X)\right)}{\mu(V)},
 $$
 where we denote
 $$
 \Ball_R(\X):=\{y\in V\mid \exists\, x\in \X, \ \dist(x,y)\le R\}.
 $$ 
 In other words, $\D_R(\X)$ is  the fraction of the measure space $V$ covered by closed radius-$R$ balls around the points of~$\X$. The problem of determining or estimating the \emph{$R$-packing density} $\D_R(V)$ of a given metric space $V$ (that is, the supremum of the densities of $R$-packings in $V$) is the archetypical problem of coding theory which also has a large number of applications to other fields.
 
One important case is when $V=\I R^m$, endowed with the Euclidean distance
and the Lebesgue (uniform) measure.
Since the measure is not finite here, one defines the \emph{$R$-packing density} $\D(\I R^m)$ as the limit of the packing densities of growing cubes that exhaust $\I R^m$:
 \beq{eq:limit}
 \D(\I R^m):= \lim_{n\to\infty}\D_R\left(\,[-n,n]^{m}\,\right).
 %\D_R(([-n,n]^m,\dist_2,\lambda).
\eeq
 It is easy to see that the limit exists and, by the scaling properties of the Lebesgue measure, does not depend on the radius~$R$. 
 %We refer the reader to the book by Conway and Sloane~\cite{ConwaySloane:splg} for an extensive treatment of this problem (prior to 1999).

It is clear that the packing density of the real line $\D(\I R^1)=1$.
The packing density of $\I R^m$ for $m=2,3,8,24$ was determined by respectively Thue~\cite{Thue92}, Hales~\cite{Hales05}, Viazovska~\cite{Viazovska17}, and  Cohn, Kumar, Miller, Radchenko and Viazovska~\cite{CohnKumarMillerRadchenkoViazovska17}.%
\hide{ Thue~\cite{Thue92}, in 1892, proved that $\D(\I R^2)=\pi/\sqrt{12}=0.9068\dots$\,, which is achieved by the hexagonal lattice; Hales~\cite{Hales05} proved in 2005 that $\D(\I R^3)=\pi/\sqrt{18}=0.7404...$\,. The packing density for $m=8$ and $m=24$ was determined by respectively Viazovska~\cite{Viazovska17} and Cohn, Kumar, Miller, Radchenko and Viazovska~\cite{CohnKumarMillerRadchenkoViazovska17}.}
The value of $\D(\I R^m)$ is still unknown for any other $m$, although the bounds for $m=4$ (from~\cite{KorkineZolotarev72,CohnLaatSalmon22}) are rather close to each other. We refer the reader to Cohn~\cite{Cohn24SpherePackingTable} for the table of the known bounds for $m\le 48$.

When $m\to\infty$, the best known upper and lower bounds on $\D(\I R^m)$ are exponentially far apart. The best known upper bound $\D(\I R^m)\leq 2^{-(0.5990...+o(1))\, m}$ is due to  Kabatjanski\u{\i} and Leven\v{s}te\u{\i}n~\cite{KabatjanskiiLevenstein78}, see also~\cite{CohnElkies03,CohnZhao14}. 
%As shown by Cohn and Zhao~\cite{CohnZhao14}, the more direct approach of Cohn and Elkies~\cite{CohnElkies03} gives at least as strong upper bound on $\D(\I R^d)$ as that of~\cite{KabatjanskiiLevenstein78}.

The trivial \emph{covering lower bound} $\D(\I R^m)\geq 2^{-m}$ (take a maximal packing and observe that balls of doubled radius cover the whole space) was improved by a factor of $\Omega(m)$ by Rogers~\cite{Rogers47}. After a sequence of improvements in the constant factor~\cite{DavenportRogers47,Ball92,Vance11,Venkatesh13}, Venkatesh~\cite{Venkatesh13} was the first to beat the $\Omega(\frac{m}{2^{m}})$ lower bound by showing that  $\D(\I R^m)=\Omega(\frac{m\, \ln\ln m}{2^{m}})$ for some infinite (sparse) sequence of dimensions~$m$. In a recent breakthrough, Campos, Jenssen, Michelen and Sahasrabudhe \cite{CamposJenssenMichelenSahasrabudhe} proved that $\D(\I R^m)\ge (\frac12-o(1))\, \frac{m\, \ln m}{2^{m}}$. This was done by discretising the problem and then applying their new result (which is stated as Theorem~\ref{thm: indep} here) on the existence of a large independent set in graphs with sparse neighbourhoods. 

The aim of this paper is to observe that the method from~\cite{CamposJenssenMichelenSahasrabudhe} also applies to ball packing in hyperbolic spaces.  Let $\I H^m$ denote the $m$-dimensional hyperbolic space, equipped with the standard metric $\dist_m$
% making it of constant curvature $-1$ 
and the standard (isometry invariant) measure $\HMeasure$.
% coming from the element of distance. 
(See Section~\ref{se:Hm} for all formal definitions.) 

% needed here, and e.g.\ Bridson and Haefliger~\cite[Section~2]{BridsonHaefliger:msnpc} or Ratcliffe~\cite{Ratcliffe19fhm} for detailed treatment of hyperbolic spaces (including their aspects relevant to this paper).

The study of packings in $\I H^m$ was held back for a long while by the absence of a good definition of the maximum packing density, see e.g.\ the discussion in \cite[Pages 831--834]{FejestothKuperberg93}. 
Analogues of~\eqref{eq:limit} lack many desired properties because growing polygons or balls in $\I H^m$ have most of their mass near their boundary. 
So some ``local'' approaches were considered; see~\cite[Chapter 11]{FejestothFejestothKuperberg23lapss} for a survey  of known results from this point of view. 
However, even here one has to be careful in view of the following striking example of B\"or\"oczky~\cite{Boroczky74} (also described on Page 833 of \cite{FejestothKuperberg93}). Namely, B\"or\"oczky showed that
% for some values of $R$, 
there is an $R$-packing $\X\subseteq\I H^2$ in the hyperbolic plane  and two tilings $\C T_1$ and $\C T_2$ of $\I H^2$ by single polygonal tiles $T_1$ and $T_2$ respectively such that each tile of $\C T_i$ contains one radius-$R$ disk centred at a point of $\X$ and is disjoint from all other such disks; at the same time, $T_1$ and $T_2$ have different areas. Thus if we measure the ``density'' of this packing $\X$ using the fraction occupied by the corresponding disks within each tile of $\C T_1$ (resp.\ $\C T_2$), then we get two different numbers. 
%This rules out various ``local'' definitions of density of a single packing.

A  solution was proposed by Bowen and Radin~\cite[Definition~5]{BowenRadin03} who presented a new notion of density where, roughly speaking, one considers probability distributions on $R$-packings in $\I H^m$ that are invariant under isometries and maximises the probability that a fixed point is covered by a ball. This parameter has many nice properties and is the one that will be used in this paper as the definition of the \emph{$R$-packing density} of $\I H^m$ and thus denoted by  $\D_R(\I H^m)$. We refer the reader to Section~\ref{se:D} for the formal definition of $\D_R(\I H^m)$ and to \cite{Radin04,BowenRadin04} for further discussions and motivation behind it.

For dimension $m=2$ and countably many radii $R$, Bowen and Radin~\cite[Theorem~2]{BowenRadin03} were able to determine $\D_R(\I H^2)$. Namely, for every integer $n> 6$, if we tile $\I H^2$ by (equilateral) triangles having all three angles $2\pi/n$ and take a ``uniform shift" of the set of the triangles' vertices then this distribution attains $\D_R(\I H^2)$, where $R=R_n$ is the half of the side length of the triangles.
As far as we know, these are the only pairs $(m,R)$ with $m\ge 2$ for which the exact value of $\D_R(\I H^m)$ is known.

Let us discuss the known upper bounds on $\D_R(\I H^m)$ for $m\to\infty$.
Before Bowen and Radin's paper~\cite{BowenRadin03}, Fejes T\'oth~\cite{Fejestoth43}  (for $m=2$), B\"or\"oczky and Florian~\cite{BoroczkyFlorian64} (for $m=3$) and B\"or\"oczky~\cite{Boroczky78} (any $m$) proved that for every $R$-packing $\X\subseteq\I H^m$ and any $x\in \X$ the fraction of volume occupied by the \emph{$R$-ball} $\Ball_R(x):=\Ball_R(\{x\})$ around $x$
%$$
%\Ball_R(x):=\{y\in\I H^m\mid \dist_m(x,y)\le R\}
%$$
inside the Dirichlet-Voronoi cell of $x\in\X$ is at most $\delta_m(2R)$, which is defined to be the fraction of volume of a regular simplex of side length $2R$ occupied by (touching) balls of radius $R$ at its vertices. It follows from these results (see Lemma~\ref{lm:Boroczky} here) that $\D_R(\I H^m)\le \delta_m(2R)$.
More recently, Cohn and Zhao~\cite[Theorem 4.1]{CohnZhao14} showed that, independently of $R>0$,
the Bowen--Radin density satisfies
\beq{eq:CohnZhao}
	\D_R(\I H^m)\le
	\inf_{\pi/3\le \theta\le \pi} \sin^{m-1}(\theta/2) A(m,\theta),
\eeq
 where  $A(m,\theta)$ be the maximum number of unit vectors in the Euclidean space $\I R^m$ such that the scalar product of every two is at most $\cos\theta$; in other words, $A(m,\theta)$ is the maximum size of a spherical code with minimum angle~$\theta$. In particular, the method of Kabatjanski\u{\i} and Leven\v{s}te\u{\i}n~\cite{KabatjanskiiLevenstein78} (that also applies to spherical codes) gives  that
 \beq{eq:CohnZhao2} 
  \D_R(\I H^m)\le 2^{-(0.5990...+o(1))\,m},
  \eeq
 see \cite[Corollary 4.2]{CohnZhao14}. For large $m$, this bound is smaller than $\delta_m(0):=\lim_{R\to0} \delta_m(R)$, which can be shown to be the volume ratio coming from $m+1$ touching unit balls in $\I R^m$. Note that Marshall~\cite[Theorem 2]{Marshall99} proved that, for all large $m$, the function $\delta_m(R)$ is strictly increasing in~$R$ (and thus  $\delta_m(0)\le \delta_m(R)$ for any $R>0$ for such $m$).
 \hide{Also, it is not hard to show (see e.g. \cite[Page 194]{Kellerhals98}) that, for every $m$, the limit
  $$
   \delta_m(\infty):=\lim_{R\to\infty} \delta_m(R)
   $$
   exists. Kellerhals \cite[Corollary 3.4]{Kellerhals98} gives an explicit formula for $\delta_m(\infty)$ but, even in this special case, the formula is rather complicated. 
   The closed-form  values $\delta_m(\infty)$ for $m\le 9$ can be found in~\cite[Page 55]{KozmaSzirmai23}. 
   This formula, with some obvious asymptotic simplifications,  reads:
 	$$
 	\delta_n(\infty)=(1+o(1)) \frac{\sqrt n}{2^{n-1}}\cdot
 	\frac
 	{\prod_{k=2}^{n-1} \left(\frac{k-1}{k+1}\right)^{\frac{n-k}2}}
 	{\sum_{k=0}^{\infty}\frac{\beta(\beta+1)\dots(\beta+k-1)}{(n+2k)!} A_{n,k}},
 	$$
 	where $\beta:=(n+1)/2$ and 
 	$$
 	A_{n,k}:=\sum_{i_0+\dots+i_n=k\atop i_j\ge 0} \frac{(2i_0)!\dots(2i_k)!}{i_0!\dots i_k!}.
 	$$
 }

The covering principle (that if $\X$ is a maximal $R$-packing then $2R$-balls around points of $X$ cover the whole space) translates with some work (see Lemma~\ref{lm:TrivialLower}) into the basic lower bound
 \beq{eq:TrivialL}
 \D_R(\I H^m)\ge 	L_R(\I H^m):=\frac{\HMeasure(\Ball_{R})}{\HMeasure(\Ball_{2R})},
 \eeq
where $\HMeasure(\Ball_r)$ denotes the volume of some (equivalently, any) $r$-ball in $\I H^m$.
 
 There are various constructions of packings in small dimensions (see e.g.\ \cite[Chapter 11]{FejestothFejestothKuperberg23lapss} for references) and it is plausible that their ``locally'' measured densities translate into lower bounds on the Bowen-Radin density. However, we are not aware of any improvements to~\eqref{eq:TrivialL} for large $m$ apart that, as far as we can see, the method of Jenssen, Joos and Perkins~\cite{JenssenJoosPerkins18,JenssenJoosPerkins19}, which is based on the hard-sphere model of statistical physics, can be applied to improve the bound in~\eqref{eq:TrivialL} by factor $\Omega(m)$. Anyway, this is superseded by the main result of this paper as follows.
 
 \begin{theorem}\label{th:main} For every $\eps>0$ there is $m_0$ such that for any $m\ge m_0$ it holds for every $R\in (0,\infty)$ that
  \[
 %\beq{eq:main}
 \D_R(\I H^m) \geq (1-\eps)\, m\ln\left(\sqrt{m} \cosh(2R)\right) \frac{\vol(\Ball_R)}{\vol(\Ball_{2R})},
 \]
 where $\ln$ is the natural logarithm and $\cosh$ is the hyperbolic cosine.
 \end{theorem}
 
Note that $\ln(\cosh (2R))$ is at least $2R-1$, so Theorem~\ref{th:main} improves the bound in~\eqref{eq:TrivialL} by factor at least
 $\Omega(m(R+\ln m))$. We prove Theorem~\ref{th:main} by reducing the lower bound problem to finding a packing in some finite-volume space (namely, the quotient of $\I H^m$ by a large-girth lattice) and then (as it was done in \cite{CamposJenssenMichelenSahasrabudhe}) discretising the problem by taking a Poisson point process. 
 
\medskip \noindent\textbf{Organisation of the paper.} We recall some notions related to $\I H^m$ in Section~\ref{se:Hm} and define the Bowen-Radin density in Section~\ref{se:D}. We describe the method from~\cite{BowenRadin03} of lower bounding $\D_R(\I H^m)$  via $R$-packings in some finite-volume space $M$ in Section~\ref{se:M}. Various estimates are collected in Section~\ref{se:Estimates}. Finally, Theorem~\ref{th:main} is proved in Section~\ref{se:Proof}.
 Since $m$ and $R$ will be reserved, respectively, for the dimension of the studied hyperbolic space and the packing radius, we may omit them from our notation if the meaning is clear.

\section{Hyperbolic space}\label{se:Hm}
 	
This section contains just a bare minimum of material sufficient to formally define the hyperbolic space $\I H^m$ and some needed related notions. For a detailed introduction to hyperbolic spaces see 
 	e.g.\ Bridson and Haefliger~\cite[Section~2]{BridsonHaefliger:msnpc} or Ratcliffe~\cite{Ratcliffe19fhm}. 
 	
One representation
 	of $\I H^m$ (for more details see e.g.\ \Ra{Chapter 3}) is to take the bilinear form 
 	\begin{equation}\label{eq:form}
 		\langle u, v\rangle_{m,1}:=-u_{m+1}v_{m+1}+\sum_{i=1}^m u_iv_i,\quad u,v\in\I R^{m+1},
 	\end{equation} 
 	identify $\I H^m$ with upper sheet of the hyperboloid  
 	$$
 	\C H:=\{u\in \I R^{m+1}\mid \langle u, u\rangle_{m,1}=-1\}
 	$$ (namely, the sheet where $u_{m+1}>0$), and
 	define the metric $\dist_m$ by $\cosh \dist_m(u,v)=-\langle u, v\rangle_{m,1}$ for $u,v\in \I H^m$.  
 	%Then all geodesic lines are precise intersections
 	%of $\I H^m$ with 2-dimensional vector subspaces of $\I R^d$ 
 	%(\cite[Corollary~2.8]{BridsonHaefliger:msnpc}). 
 	The group of isometries of $\I H^m$
 	can be identified with $O(m,1)_0$, the group of $(m+1)\times (m+1)$-matrices $A$ which leave
 	the bilinear form in~\eqref{eq:form} invariant
 	% (that is, $A^TJA=J$, where $J$ is the diagonal matrix with  the first $m$ diagonal entries $+1$ and the last one $-1$)
 	and do not swap the two sheets of $\C H$; see~\cite[Theorem~2.24]{BridsonHaefliger:msnpc} or~\cite[Theorem~3.2.3]{Ratcliffe19fhm}.

 	Let $\CG$ be the group  of orientation-preserving isometries of $\I H^m$; it corresponds to the 
 	subgroup $\SO(m,1)_0$ of index 2 in $O(m,1)_0$,  which 
 	consists of 
 	matrices with determinant 1. It is a topological group with the topology inherited from $O(m,1)_0\subseteq \I R^{(m+1)^2}$.
 	
\hide{ Recall that the \emph{$r$-ball} around $x\in\I H^m$ is 
 	$$
 	\Ball_r(x):=\{y\in\I H^m\mid \dist_m(x,y)\le r\}.
 	$$ 
 	%Also, we denote $\Ball_r:=\Ball_r(\C O)$.
 }
 
	We fix one point of $\I H^m$, say $(0,\dots,0,1)\in \I R^{m+1}$ in the above representation, which we will denote by $\CO$ and refer to as the \emph{origin}. 
	
 	The space $\I H^m$ is equipped with a Borel isometry-invariant measure $\HMeasure$ whose push-forward
 	under the projection on the first $m$ coordinates of $\I R^{m+1}$ has density 
 	$(1+(x_1^2+\ldots+x_m^2))^{-1/2}$
% 	\beq{eq:HnDensity}
% 	\rho(x_1,\ldots,x_m):=(1+(x_1^2+\ldots+x_m^2))^{-1/2}
% 	\eeq
 	with respect to the Lebesgue measure on~$\I R^m$. 
\hide{The real lie algebra $\mathfrak{so}(n,1)$ of the group $\CG$ is simple (e.g.\ by the full characterization that can be found in~\cite[Theorem 6.105]{Knapp:lgbi}) so the group is unimodular (e.g. by~\cite[Corollary 5.5.5]{Faraut08alg}). More directly (but with some calculations required) one can consider the explicit description of $\mathfrak{so}(n,1)$ as algebra of matrices $M$ with $M^TI_{n,1}+I_{n,1}M=0$ and check that the adjoint action of $\GG$ on it by conjugations has determinant $1$.} 
The group $\CG$ is unimodular and locally compact, so there is a Haar measure $\mu_{\CG}$ on $\CG$ (which is both right- and left-invariant).  It is a standard result (see e.g.\ \Ra{Lemma 4 of Section 11.6}) that, by scaling $\mu_{\CG}$, we can assume that the projection map $\pi_{\CO}:(\CG,\mu_{\CG})\to (\I H^m,\HMeasure)$, where  $\gamma\mapsto \gamma.\C O$, is measure-preserving.
 	
 	We consider the natural (left) action $\CG\actson\I H^m$. For $\gamma\in \CG$ and $A\subseteq \I H^m$, we denote 
 	%$\gamma.x:=a(\gamma,x)$ and 
 	$\gamma.A:=\{\gamma.x\mid x\in A\}$. For $x\in \I H^m$, let 
 	$
 	\Sigma_x:=\{\gamma\in\CG\mid \gamma.x=x\}
 	$
 	be the \emph{stabiliser} of~$x$. In particular, we have the natural homeomorphism of topological spaces $\CG/\Sigma_{\C O}\cong\I H^m$.

	\section{Definition of the Bowen--Radin density}\label{se:D}
	
In this section, we give the definition of the packing density $\D_R(\I H^m)$ introduced by Bowen and Radin~\cite{BowenRadin03}. We also state some results from~\cite{BowenRadin03}, occasionally providing more details (usually when these were implicitly assumed but not stated in~\cite{BowenRadin03}).

Let $R>0$ be a positive real. Let $S_R$ consist of $R$-packings $\X\subseteq \I H^m$ that are \emph{relatively dense}, that is, such that, for every $x\in\I H^m$, the (closed radius-$2R$) ball $\Ball_{2R}(x)$ around $x$ contains at least one element of~$\X$. 
This is slightly weaker than the notion of a maximal packing (namely, a relatively dense packing can have  the distance from some $x$ to $\X$ exactly $2R$ and thus not be a maximal one). 
	
Bowen and Radin~\cite{BowenRadin03} considered the following metric $\dist_R$ on $S_R$:
	\beq{eq:dR}
	\dist_R(\X,\Y):=\sup_{r\in[1,\infty)} \frac1r\, h\big(\Ball_r(\CO)\cap \X,\Ball_r(\CO)\cap \Y\big),\quad \X,\Y\in
	S_R,
	\eeq
	where 
	$$
	h(A,B):=\max\left\{ \sup_{a\in A}\inf_{b\in B} \dist_m(a,b),\ \sup_{b\in B}\inf_{a\in A}\dist_m(a,b)\right\}, \quad A,B\subseteq \I H^m,
	$$
	is the usual \emph{Hausdorff distance} on the subsets of $\I H^m$. As noted in~\BR{Page 25}, the metric space $(S_R,\dist_R)$ is compact. Let $\C M(R)$ be the set of Borel probability measures on $(S_R,d_R)$. 
	
The group $\CG$ naturally acts on $S_R$ and the corresponding map $\CG\times S_R\to S_R$ is continuous.
	A measure $\mu\in \C M(R)$ is called \emph{$\CG$-invariant} if for every $\gamma\in \CG$ and every Borel set $E\subseteq S_R$ we have $\mu(\gamma.E)=\mu(E)$. Let $\C M_I(R)$ consist of all $\CG$-invariant measures in $\C M(R)$. Let $\C M_I^{\mathbf{e}}(R)$ consist of those  measures $\mu\in \C M_I(R)$ that are \emph{ergodic}, that is, are extreme points of the convex set $\C M_I(R)$.
	
	For $p\in\I H^m$, define 
	 $$
	 F_p(\X):=\left\{\begin{array}{ll} 1,& \mbox{if $\dist_m(x,\X)\le R$,}\\
	 0, & \mbox{otherwise},\end{array}\right.
	 \quad \text{ for } \X\in S_R,
	 $$
	  to be the indicator function that $R$-balls around $\X$ cover~$p$.

	\begin{lemma}\label{lm:FpMeas} 
		For every $p\in\I H^m$, the function $F_p$ is a Borel function on $(S_R,\dist_R)$.
	\end{lemma}
	\bpf Let us show that $F_p^{-1}(1)=\{\X\in S_R\mid F_p(\X)=1\}$ is a closed subset of~$(S_R,\dist_R)$. Take any sequence $(\X_n)_{n=1}^\infty$ in $F_p^{-1}(1)$ convergent to some $\X\in S_R$. Each $\X_n$ has a point $x_n$ in $\Ball_R(p)$. By the compactness of the closed ball $\Ball_R(p)$, we can pass to a subsequence of $n$ so that $x_n$ converges to some $x\in\Ball_R(p)$ as $n\to\infty$. Since $\X$ has at most one point in $\Ball_{R/2}(x)$, it follows from the definition of the Hausdorff distance by taking any $r>\dist_m(x,\CO)$ in~\eqref{eq:dR} that $x\in \X$. Thus $F_p(\X)=1$ and the set $F_p^{-1}(1)$ is closed (and thus Borel). 
	
	Since the function $F_p$ is $\{0,1\}$-valued and the pre-image $F_p^{-1}(1)$ is closed, the function $F_p$ is Borel.\epf
	
	Thus for any $\mu\in \C M_I(R)$, we can define the \emph{average density} 
	\beq{eq:Dmu}
	D(\mu):=\int_{S_R} F_{\C O}(\X)\dd\mu(\X).
	\eeq
	By the $\CG$-invariance of $\mu$ and the transitivity of $\CG\actson\I H^m$, we could have taken any other point of $\I H^m$ instead of the origin in this definition, without changing the value.
	
Following Bowen and Radin~\BR{Definition~5}, we define the \emph{$R$-packing density of $\I H^m$} as
	\beq{eq:DR}
	\D_R(\I H^m):=\sup_{\mu\in\C M_I^{\mathbf{e}}(R)} D(\mu).
	\eeq
	We refer the reader to \cite{BowenRadin03} and also to~\cite{BowenRadin04,Radin04} for discussion and further properties of this parameter. For example, \BR{Theorem~1} implies that the supremum in~\eqref{eq:DR} is attained by some measure $\mu\in\C M_I^{\mathbf{e}}(R)$. 
	
Let us show that the definition in in~\eqref{eq:DR} is not affected if we take the supremum over invariant (not necessarily ergodic) measures.	

\begin{lemma}\label{lm:RemoveErgodic}
 For every $R\in (0,\infty)$ and $m\in\I N$, we have that $\D_R(\I H^m)=\sup_{\mu\in\C M_I(R)} D(\mu)$. 
\end{lemma}

\bpf
Trivially, $\D_R(\I H^m)\le \sup_{\mu\in\C M_I(R)} D(\mu)$, so let us show the converse
inequality. 

The ergodic invariant decomposition theorem of Varadarajan~\cite[Theorem 4.2]{Varadarajan63} (see also Farrell~\cite[Theorem 5]{Farrell62} for a similar result) applies to Borel actions of locally compact groups on standard Borel spaces. When applied to
% the measure $\mu$ invariant under 
the action $\C G\actson S_R$, it gives a Borel map $\beta:S_R\to\C M_I^{\mathbf{e}}(R)$ which is constant on every orbit of $\CG$ such that  every measure $\nu$ in the image of $\beta$ assigns value $1$ to the pre-image $\beta^{-1}(\nu)\subseteq S_R$ and every $\mu\in\C M_I(R)$ satisfies
 \beq{eq:Varadarajan}
 \mu(A)=\int_{S_R} \beta(\X)(A)\dd\mu(\X),\quad\mbox{for every Borel $A\subseteq S_R$}.
 \eeq

Take any $\mu\in \C M_I(R)$ and let $A$ be $\{\X\in S_R\mid \C O\in\Ball_R(\X)\}=F_{\CO}^{-1}(1)$. By Lemma~\ref{lm:FpMeas}, $A$ is a Borel subset of $S_R$.  Note that $\mu(A)=\D(\mu)$. Thus, by~\eqref{eq:Varadarajan} applied to $\mu$ and $A$, $\D(\mu)$ is some average of $\D(\nu)$ over ergodic invariant measures $\nu=\beta({\X})$ for $\X\in S_R$. Thus there is $\nu\in  \C M_I^{\mathbf{e}}(R)$ with $\D(\nu)\ge\D(\mu)$ and we have $\D_R(\I H^m)\ge \D(\mu)$. 
Since $\mu\in \C M_I(R)$ was arbitrary, the desired inequality follows.
\epf

%It can be shown that the definition in~\eqref{eq:DR} does not change if we remove the requirement that $\mu$ is ergodic (that is, take the supremum over the larger set $\C M_I(R)$). 
Also, it is plausible that one can replace $S_R$ by the space of arbitrary (not necessarily relatively dense) $R$-packings by completing a random $R$-packing into a relatively dense one in a measurable and invariant way. However, this would require some extra work (and the definition of the distance in~\eqref{eq:dR} would need some tweaking) so we stay with the definitions from~\cite{BowenRadin03}.

At this point, let us observe
% (since this was not explicitly stated in~\cite{BowenRadin03}) 
that the classical ``local'' upper bounds from~\cite{Fejestoth43,BoroczkyFlorian64,Boroczky78} also apply to $\D_R(\I H^m)$.

\begin{lemma}\label{lm:Boroczky} For any real $R\ge 0$ and integer $m$, we have $\D_R(\I H^m)\le \delta_m(2R)$.\end{lemma}
	
\bpf Take any $\mu\in\C M_I^{\mathbf{e}}(R)$. Bowen and Radin~\BR{Proposition 3} showed that $D(\mu)$ can be computed as the expectation over a random packing $\X$ distributed according to $\mu$ of the ratio of volume of $\Ball_{R}(\CO)$ to the volume of the Dirichlet-Voronoi cell of $\X$ containing the origin $\C O\in\I H^m$ in its interior. (Note that, by the invariance of $\mu$, the probability that $\C O$ lies on the boundary of a Dirichlet-Voronoi cell of $\X$ is zero.) 

%By the result of B\"or\"oczky~\cite{Boroczky78}  
By the results from~\cite{Fejestoth43,BoroczkyFlorian64,Boroczky78}
mentioned in the Introduction, the integrated
function is always at most $\delta_m(2R)$, so its average is at most $\delta_m(2R)$, as desired.\epf

\section{Lower bounds via a finite-volume space}\label{se:M}

Here we describe the method, implicit in~\cite{BowenRadin03}, of proving lower bounds on $\D_R(\I H^m)$ coming from ball packings in some finite-volume space. We use (as in~\cite{KuperbergKuperbergFejestoth98}) the following definition of a lattice which is stronger than the standard one.
% assumed in~\cite{BowenRadin03} (so all cited positive results from~\cite{BowenRadin03} involving lattices remain true).
Namely, let us call a subgroup
	$\CL\subseteq \CG$  a \emph{(uniform torsion-free) lattice} if (i) there is $g>0$ such that for every $x\in\I H^m$ and every non-identity $\gamma\in\CL$ we have $\dist_m(x,\gamma.x)\ge g$; (ii) the quotient space $M:=\I H^m/\CL$ is compact. Denote the supremum of $g$ that work in (i) as $g(\CL)$ and call it the \emph{girth} of $\CL$. The topological space $M$ comes with the natural metric 
	$$
	\dist_M(x,y):=\dist_m(\pi_M^{-1}(x),\pi_M^{-1}(y)),\quad x,y\in M,
	$$ 
	where $\pi_M:\I H^m\to M$ is the projection. Note that $\pi_M$  is a \emph{local isometry}, that is, there is $r>0$ such that for every $x\in \I H^m$ the map $\pi_M$ gives an isometry between $\Ball_r(x)$ and the ball of radius $r$ around $\pi_M(x)$ in $\dist_M$. 
	
	Let us show that, in fact, we can take the local isometry radius $r$ to be $\frac14\, g(\CL)$.
	Take any $x\in\I H^m$. To show the surjectivity on radius-$r$ balls, take any $w\in M$ with $\dist_M(w,\pi_M(x))\le r$. Since $\pi_M^{-1}(\pi_M(x))=\C L.x$ and $\C L$ acts by isometries, we have
	$
	%\dist_m(\C L.x,\pi_M^{-1}(w))=
	\dist(\pi_M^{-1}(w),x)=\dist_m(\pi_M^{-1}(w),\C L.x)\le r$. Since $\pi_M^{-1}(w)$ is closed and, say, $\Ball_{2r}(x)$ is compact, we have that $\pi_M^{-1}(w)\cap \Ball_r(x)$ is non-empty. Thus $\pi_M$ is a surjection of $\Ball_{r}(x)$ onto the ball of radius $r$ around $\pi_M(x)$. To show that $\pi_M$ preserves distances on $\Ball_r(x)$, take any $y,z\in\Ball_r(x)$. Clearly, $d_m(y,z)\ge d_M(\pi_M(y),\pi_M(z))$. For the converse inequality, we have to show that $d_m(y',z')\ge d_m(y,z)$ for any $y',z'\in\I H^m$ with $\pi_M(y')=\pi_M(y)$ and $\pi_M(z')=\pi_M(z)$. By applying an element of $\C L$
	% (which is an isometry of $\I H^m$) 
	to $y'$ and $z'$, we can assume that $z'=z$. If $y=y'$ then there is nothing to do; otherwise $d_m(y,y')\ge g(\CL)$ and thus $d_m(y',z)\ge d_m(y',y)-d_m(y,z)\ge g-2r \ge 2r\ge d_m(y,z)$, as required.
	
	 Since $\I H^m$ has constant curvature $-1$ the same applies to $M$; so $M$ is a hyperbolic manifold (compact, without boundary). 
	Since the  measure $\HMeasure$ on $\I H^m$ can be defined via the element of length, see e.g.\ \Ra{Section 3.4}, this definition carries over to $M$ giving that
	% $\pi_M$	is \emph{locally measure preserving} (that is, 
	there is $r$ (which we can take again to be $\frac14\,g(\CL)$) such that for every $x\in\I H^m$ the restriction of $\pi_M$ to $\Ball_r(x)$ is measure-preserving. By the compactness of $M$, it can be covered by finitely many balls of radius $\frac14\,g(\CL)>0$. Each of these balls has finite volume (since finite-radius balls in $\I H^m$ have finite volume) so $\mu_M(M)<\infty$.

	\begin{lemma}\label{lm:CInverse} Let a lattice $\CL\subseteq \CG$ have girth $g\geq 8R$, define $M:=\I H^m/\CL$, and let $\pi_M:\I H^m\to M$ be the projection. Let $R>0$ and let $\Y$  be an $R$-packing in $(M,\dist_M)$. 
		Define $\X:=\pi_M^{-1}(\Y)\subseteq \I H^m$. Then the following statements hold.
		
		\begin{enumerate}
			\item\label{it:CInverse1} The set $\X$ is an $R$-packing in $(\I H^m,\dist_m)$.
			%\item\label{it:CInverse2} If the packing $\X$ is maximal then $\X'$ is maximal.
			\item\label{it:CInverse3} If the packing $\Y$ is relatively dense in $(M,\dist_M)$ then $\X$ is relatively dense in $(\I H^m,\dist_m)$.
	\end{enumerate} \end{lemma}
	\bpf %Suppose that $\X$ is a maximal packing.
	Let us show the first claim that $\X$ is a packing. Take any distinct $x,x'\in\X$ and let $y:=\pi_M(x)$ and $y':=\pi_M(x')$. If $y\not=y'$, then $\dist_M(y,y')\ge 2R$ and, by the definition of $\dist_M$, we have $\dist_m(x,x')\ge 2R$. If $y=y'$ then, since $\pi_M$ maps $\Ball_{g/4}(x)$ isometrically to the ball of radius $g/4\geq 2R$ around $y$ in $M$ and this ball contains $\pi_M(x')=y$, we have $\dist_m(x,x')> 2R$, giving that $X$ is an $R$-packing, as desired. 

Let us show  that the packing $\X\subseteq \I H^m$ is relatively dense if $\Y$ is. Take any $x\in\I H^m$. Since the packing $\Y$ is relatively dense, we have that $\dist_M(y,\Y)\le 2R$, where $y:=\pi_M(x)$. 
 %Since finite-radius balls in $\I H^m$ are compact, there is $y'\in\Y$ with $\dist_M(y,y')\le 2R$. 
As $\C L$ acts by isometries and $\X$ is invariant under $\C L$, we have
	$$
	\dist_m(x,\X)= 
	\dist_m(\C L.x,\X)=
	\dist_m(\pi_M^{-1}(y),\pi_M^{-1}(\Y)) 
	=\dist_M(y,\Y)
	%\le \dist_M(y,y')
	\le 2R.
	$$
	By the compactness of, say, $\Ball_{3R}(x)$, we have that $\Ball_{2R}(x)\cap X\not=\emptyset$.
	Since $x\in\I H^m$ was arbitrary, the packing $\X$ is indeed relatively dense.
\epf

\begin{lemma}\label{lm:Transfer} If $\CL\subseteq\CG$ is a lattice of girth at least $8R$ and $M=\I H^m/\C L$ then
 \beq{eq:Transfer}
 \D_R(\I H^m)\ge \D_R(M).
 \eeq
 \end{lemma}
\bpf  
 Take any relatively dense $R$-packing $\Y$ in $M$ and let $\X:=\pi_M^{-1}(\Y)$. Then $\X$ is a relatively dense $R$-packing in $\I H^m$ 
by Lemma~\ref{lm:CInverse}. Moverover, it is invariant under the action of the lattice $\CL\actson\I H^m$.

Let $\Gamma_{\X}:=\{\gamma\in\CG\mid \gamma.\X=\X\}$ be the subgroup of $\CG$ that fixes~$\X$. Clearly, $\C L\subseteq \Gamma_{\X}$. 

Define the probability measure $\mu_X$ on $S_R$ by 
$$
 \mu_{\X}(E):=\mu_{\CG}(\{\gamma\in\CG\mid \gamma.\X\in E\}),\quad \mbox{for Borel $E\subseteq S_R$}.
 $$
 Informally speaking, we take the translate of the given $R$-packing $\X$ by a random element of $\CG$.
 
 This measure is $\CG$-invariant. Indeed, for every $\gamma'\in\CG$ and $E\subseteq S_R$, we have by the invariance of $\mu_{\CG}$ that
 \begin{eqnarray*}
 \mu_{\X}(\gamma'.E)&=&\mu_{\CG}(\{\gamma\in\CG\mid \gamma.\X\in \gamma'.E\})\\
 &=&\mu_{\CG}(\gamma'.\{\gamma''\in\CG\mid \gamma''.\X\in E\})\\
 &=&\mu_{\CG}(\{\gamma''\in\CG\mid \gamma''.\X\in E\})\ =\ \mu_{\X}(E).
 \end{eqnarray*}

%By $\mu_{\X}\in\C M_I(R)$ and 
By Lemma~\ref{lm:RemoveErgodic} we have that $\D_R(\I H^m)\ge \D(\mu_{\X})$.

By  \BR{Proposition 1}, the density $D(\mu_{\X})$ is equal to the relative volume taken by the $R$-balls around $\X$ inside any fundamental domain of $\Gamma_{\X}\actson\I H^m$. The lattice $\CL$, as a subgroup of $\Gamma_{\X}$, has finite index (which can be upper bounded by $k!$ where $k$ is the maximum size of an $R$-packing in $M$, where $k$ in turn  can be upper bounded by $\mu_M(M)/\HMeasure(\Ball_R)<\infty$). Thus a fundamental domain $\C F\subseteq \I H^m$ of $\CL\actson \I H^m$ can be obtained by taking the union of finitely many (pairwise disjoint) translates of a fundamental domain of $\Gamma_{\X}\actson\I H^m$ (one per each coset of $\Gamma_{\X}/\CL$). Since each of the latter translates has the same occupied ratio (by \BR{Proposition 1}), this ratio is the same as that for their union $\C F$. Note that the restriction of $\pi_M$ to the fundamental domain $\C F$ is a measure-preserving map between $(\C F,\HMeasure)$ and $(M,\mu_M)$. By the girth assumption on $\CL$, the $R$-balls around points of the $\CL$-periodic tiling $\X\subseteq \I H^m$ occupy the same fraction of volume of $\C F$ as the $R$-balls around points of $\Y$ in~$M$. 

Since $D(\mu_{\X})$ is the radius-$R$ density of $\Y$ in $M$,
% and any packing extends to a relatively dense one, 
the lemma follows.\epf

 	\begin{lemma}\label{lm:cover} For every $m\in\I N$ and $r\in(0,\infty)$, there is a lattice $\CL$ of isometries of $\I H^m$ of girth at least $r$.
 	\end{lemma}
 	
 	\bpf This is Theorem 4.1 in~\cite{KuperbergKuperbergFejestoth98} that contains a  detailed proof. (The authors of~\cite{KuperbergKuperbergFejestoth98} wrote that this fact had been known for a long time but they could not find a suitable reference.)%
 	\hide{There exists a co-compact lattice $L$ in $\CG$, see e.g.\ \Ra{Section 12.7}. 
 		Selberg~\cite{Selberg60} proved that all lattices are virtually torsion-free. So without loss of generality, L is torsion-free. Thus all nontrivial elements of L are hyperbolic and have positive translation length. There can only be finitely many conjugacy classes of elements with translation length less than $2r$.
 		
 		All lattices are residually finite by Malcev’s Theorem~\cite{Malcev40} which states that all finitely generated linear groups are residually finite.  So there exists a finite-index subgroup $L’$ in $L$ such that all nontrivial elements of $L’$ have translation length greater than $2r$. This implies that the injectivity radius of the quotient $H^m/L’$ is at least $r$ (so the covering map $c:H^m\to M$ is an isometry on every r-ball).}
 	\epf

Recall that we defined $L_R(\I H^m):=\vol(\Ball_{R})/\vol(\Ball_{2R})$ to be the ratio of the volumes of balls of radius $R$ and $2R$ in~$\I H^m$. We can now argue that this gives a lower bound on the Bowen--Radin packing density, just to show that this natural lower bound also holds in this framework. 
 	
 	\begin{lemma}\label{lm:TrivialLower}
 	For every $m\in\I N$ and $R>0$, we have
 		$\D_R(\I H^m)\ge L_R(\I H^m)$.
 	\end{lemma}
 	
 	\bpf Take any lattice $\CL\subseteq\CG$ of girth at least $8R$, which exists by Lemma~\ref{lm:cover}. Take a maximal packing $\Y$ in $M:=\I H^m/\CL$. By the maximality of $\Y$, the $2R$-balls around points of $\Y$ cover the whole space~$M$. By the girth assumption, for any $r\le 2R$, the volume of any $r$-ball in $M$ is the same as the volume of an $r$-ball in~$\I H^m$. Thus $\D_R(M)\ge L_R(\I H^m)$. Now the lemma follows from~\eqref{eq:Transfer}.
 	\epf

	\section{Various estimates}\label{se:Estimates}
	
We will use various facts about hyperbolic functions $\cosh x:=(\me^x+\me^{-x})/2$, $\sinh x:=(\me^x-\me^{-x})/2$, $\tanh x := \frac{\sinh x}{\cosh x}$ and $\coth x:= 1/\tanh x$. First, we have $\cosh^2x=1+\sinh^2x$. The following formulas for angle doubling are easy to check directly:
	$$
	\cosh(2x)=2 \cosh^2x-1%=1+2\sinh^2x
	\quad\mbox{ and }\quad
	\sinh(2x)=2\sinh x\cosh x.
	$$
Also, we will use the monotonicity of the above hyperbolic functions on $[0,\infty)$, approximations $\tanh x,\sinh x=(1+o(1))x$ for $x\to 0$ and the following inequalities that are routine to check:
 \begin{equation}\label{eq:HyperbolicIneq}
  \tanh x< 1,\ \ 2\sinh x \le \sinh(2x),\ \ \tanh (2x)\le 2\tanh x,\quad \mbox{for $x\in [0,\infty)$.}
 \end{equation}

\hide{Also, it is routine to check that $\tanh (2x)/\tanh x$ is a decreasing function for $x\in (0,\infty)$ (namely, its derivative is $-2\tanh (2x)/\cosh (2x)<0$) approaching $2$ (resp.\ $1$) as $x\to 0$ (resp.\ $x\to\infty$) so it holds that
  \begin{equation}\label{eq:TanhRatio}
  1 < \frac{\tanh (2x)}{\tanh x}< 2, \quad\mbox{for $0<x<\infty$.}
  \end{equation}
  }
  
Furthermore, we will use the following formula (see e.g.\ \Ra{Theorem 5.3.5} for a proof).
	
	\begin{lemma}[The First Law of Cosines]
		\label{lm:1LC}
		If $\alpha,\beta,\gamma$ are the angles of a hyperbolic triangle and $a,b,c$ are the lengths of the opposite sides, then
		\beq{eq:1LC}
		\cos \gamma=\frac{\cosh a\cosh b-\cosh c}{\sinh a\sinh b}.
		\eeq
	\end{lemma}

  	\begin{lemma}\label{lm:Intersection}
		Let $x,u\in\I H^m$ be two distinct points and let $\tau:=\dist_m(x,u)>0$. Let $r>0$
		satisfy $\tau <2r$. 
		Then the intersection $\Ball_r(x)\cap \Ball_r(u)$ is contained in  (and thus has volume at most as that of) a hyperbolic ball of radius $\sigma=\sigma(\tau,r)$, where $\sigma>0$ is defined by
		\beq{eq:Intersection}
		\sinh^2(\sigma) %=\sinh^2 R-\left(\frac{\cosh R}{\sinh \tau}(\cosh \tau-1)\right)^2
  = \sinh^2(r) - \cosh^2(r) \tanh^2\left(\tau/2\right).
		% = \sqrt{\cosh^2 \tau-2\left(\frac{\cosh^2\tau-\cosh \theta}{\sinh \tau}\right)^2}.
		\eeq
	\end{lemma}
 \bpf
 Let $w$ be the point on the geodesic line through $x$ and $u$ such that $\dist_m(x,w)=\dist_m(w,u)$. We will show that $\Ball_\sigma(w)$ works in the lemma.
Let $y \in \Ball_r(x)\cap \Ball_r(u)$ be any point at distance $r_1$ from $x$ and $r_2$ from $u$ for some $r_1,r_2\leq r$. By the convexity of $\Ball_r(x)\cap \Ball_r(u)$, we can assume that $y$ is not equal to $x$ nor $u$, i.e.\ that $r_1,r_2>0$.
Consider the hyperbolic triangle with the vertices $x,y,u$ and let $\alpha$ be the angle at the vertex $u$ in this triangle. Let $\rho = \dist_m(y,w)$ be the distance between $y$ and $w$.

Applying Lemma~\ref{lm:1LC} to the hyperbolic triangle $xyu$ and the hyperbolic triangle $ywu$, we obtain 
$$ \cos \alpha = \frac{\cosh \tau \cosh r_2 - \cosh r_1}{\sinh \tau \sinh r_2} \enspace \text{ and } \enspace \cos \alpha = \frac{\cosh (\tau/2) \cosh r_2 - \cosh \rho }{ \sinh (\tau/2) \sinh r_2}.$$
As these two expressions have the equal value, using the equalities $\sinh \tau = 2\sinh (\tau/2)\cosh(\tau/2)$ and $\cosh \tau = 2\cosh^2(\tau/2) -1$, 
one can obtain 
\begin{align*}
\cosh \rho &= \cosh(\tau/2) \cos r_2 - \frac{1}{2\cosh (\tau/2)}\left(\cosh \tau \cosh r_2 -\cosh r_1\right) \\
& = \cosh(\tau/2)\cosh r_2 - \frac{(2\cosh^2(\tau/2)-1)\cosh r_2 - \cosh r_1}{2\cosh(\tau/2)} \\
& = \frac{\cosh r_2 + \cosh r_1}{2\cosh(\tau/2)} \leq \frac{\cosh r}{\cosh(\tau/2)}.
\end{align*}
The final inequality holds as $\cosh x$ is an increasing function of $x\ge 0$ and $r_1,r_2\leq r$. 
This yields 
\begin{align*}
    \sinh^2(\rho) &= \cosh^2 (\rho) -1 \leq \frac{\cosh^2(r)   }{\cosh^2(\tau/2)}-1
    = \cosh^2(r) -1 - \frac{\cosh^2(r) (\cosh^2(\tau/2) - 1 ) }{\cosh^2(\tau/2)} \\
    &= \sinh^2(r) - \cosh^2(r) \tanh^2(\tau/2).
\end{align*}
As $\sinh^2(x)$ is an increasing function on $[0,\infty)$, this shows that $\rho$ is at most $\sigma$. This proves that any point $y\in \Ball_r(x)\cap \Ball_r(u)$ has distance at most $\sigma$ from $w$, as desired.
\epf

	The volume of an $r$-ball in $\I H^m$ is given by the following formula (see e.g.\ \cite[Equation (III.4.1)]{Chavel06rgmi}): 
 \begin{equation*}
     \label{eq:VolB}
	\vol(\Ball_r)=\mathrm{vol}(S^{m-1})\int_0^r\sinh^{m-1}\eta \dd\eta,
 \end{equation*}
	where $\mathrm{vol}(S^{m-1})$ denotes the $(m-1)$-dimensional (surface) measure of the unit sphere in the Euclidean space~$\I R^m$.
	Recall that
	$$
	\mathrm{vol}(S^{m-1})=\frac{2\pi^{m/2}}{\Gamma(m/2)}=\me^{-(1/2+o_m(1)) m\ln m}.
	$$
Moreover, we have for $r\rightarrow \infty$ and $m\geq 3$ that 
\begin{align*}
    \int_0^r\sinh^{m-1}\eta \dd\eta 
&=  (1+O(\me^{-r(m-1)/2})) \int_{r/2}^r 
\left(\frac{\me^{\eta}- \me^{-\eta}}{2}\right)^{m-1} \dd\eta \\
&= (1+ O(m \me^{-r})) \int_{r/2}^r \frac{\me^{(m-1)\eta}}{2^{m-1}} \dd\eta = (1+ O(m \me^{-r}))\, \frac{e^{(m-1)r}}{(m-1)2^{m-1}}.
\end{align*}
Using these, for $r\to\infty$ and any $m\ge 3$, 
 \begin{equation}
     \label{eq:VolBr}
	\vol(\Ball_r)= (1+ O(m \me^{-r}))\frac{\me^{(m-1)r}}{(m-1)2^{m-1}}\,  \mathrm{vol}(S^{m-1}) =\me^{(m-1)r-(1/2+o_{r}(1))m\ln m}.
 \end{equation}

We will also need the following bound on the ratio between $\vol(\Ball_r)$ and $\vol(\Ball_{R})$ from \cite{CohnZhao14}*{Lemma 4.6}:

\begin{lemma}\label{lem: vol ratio}
For $0<r<R$ and the balls $\Ball_r, \Ball_R$ of radius $r,R$ in $\I H^m$, we have
  \[   \left(\frac{\sinh r}{\sinh R}\right)^m \leq \frac{\vol(\Ball_r)}{\vol(\Ball_R)} \leq  \left(\frac{\sinh r}{\sinh R}\right)^{m-1}.\]
\end{lemma}

	Let us estimate the ratio $L_R(\I H^m)=\vol(\Ball_R)/\vol(\Ball_{2R})$.
 Using the first equality in \eqref{eq:VolBr},  we have for $R$ satisfying $m \me^{-R}\rightarrow 0$ that
\begin{eqnarray*}
	L_R(\I H^m)&=& \frac{ (1+O(m\me^{-R}) )\, \me^{R(m-1)}}{ (1+O(m\me^{-2R}) )\, \me^{2R(m-1)}}
  = (1+O(m\me^{-R}) )\, \me^{-R(m-1)}.\label{eq:LEst}
	\end{eqnarray*}
\hide{ If $m \me^{- R/2} \rightarrow 0$, then $\int_{0}^{R/2} \sinh^{m-1} \eta \dd \eta \leq O(\me^{-R/2}) \int_{0}^{R} \sinh^{m-1} \eta \dd \eta$, hence
	\begin{eqnarray}
	%\frac{\int_0^R \sinh^{d-1} \eta\dd\eta}{\int_0^{2R} \sinh^{d-1} \eta\dd\eta}=
	L_R(\I H^m)&=& (1+O(\me^{-R/2}))\frac{\int_{R/2}^R \sinh^{m-1} \eta\dd\eta }{\int_{R/2}^{2R} \sinh^{m-1} \eta\dd\eta} =  (1+O(\me^{-R/2}))
 \frac{\int_0^R (\me^{(m-1)\eta}+ O(m\,\me^{(m-3)\eta})) \dd\eta}{\int_0^{2R} (\me^{(m-1)\eta}+ O(m\,\me^{(m-3)\eta}))\dd\eta} \nonumber\\
&=&	 \frac{\frac1{m-1}\me^{(m-1) R}- O(\me^{(m-3)R})}{\frac1{m-1}\me^{2(m-1) R}- O(\me^{(m-3)2R})}
	= (1+O(m\,\me^{-R/2}))\,\me^{-(m-1)R}.\label{eq:LEst}
	\end{eqnarray}}
	(If $m\me^{-R}\not=o(1)$ then we have to be more careful with the lower order terms.)
	
	It is conceivable that, for any fixed $m\ge 2$, the function $L_R(\I H^m)$ is monotone decreasing for $R\in (0,\infty)$, so it is never larger than the lower bound $2^{-m}$ from the Euclidean case (which is the limit as $R\to 0$).

Note that Lemmas~\ref{lm:1LC}--\ref{lem: vol ratio} also hold in $M=\I H^m/\CL$ as long as all involved points are are within some ball in $M$  of radius  at most $g(\CL)/4$.

 \hide{For general $R$ and $m$, assuming $m\to \infty$, we have
\beq{eq:LEst_rough}
L = \frac{C_{m,R} \sinh^{m-1} (R) }{ C_{m,2R}\sinh^{m-1}(2R) } = D_{m,R} \left(\frac{\sinh (R)}{ \sinh (2R)}\right)^{m-1},
\eeq
where $\Omega(1/m)\leq D_{m,R}\leq O(m)$.
}

 \section{Proof of Theorem~\ref{th:main}}\label{se:Proof}

Recall that Theorem~\ref{th:main} states that, for every $\eps>0$ and large $m$, we have 
\beq{eq:main} 
 \D_R(\I H^m) \geq (1-\eps) m \ln(\sqrt{m} \cosh(2R)) \frac{\vol(\Ball_{R})}{\vol(\Ball_{2R})}.
\eeq
In order to prove this, 
%the stated lower bound on $\D_R(\I H^m)$ 
we will use the following result of Campos, Jenssen, Michelen and Sahasrabudhe \cite{CamposJenssenMichelenSahasrabudhe} on independent sets in a graph $G$ with small \emph{maximum codegree} $\Delta_2(G)$, which is the maximum over distinct $x,y\in V(G)$ of the number of common neighbours of $x$ and~$y$.

\begin{theorem}\label{thm: indep}
    Let $G$ be a graph on $n$ vertices such that its maximum degree $\Delta(G)\leq \Delta$ and its maximum codegree $\Delta_2(G)\leq 2^{-7} \Delta(\ln \Delta)^{-7}$. Then the independence number $\alpha(G)$ satisfies
    $$\alpha(G)\geq (1-o(1))\,\frac{n\ln \Delta}{\Delta},$$
    where $o(1)$ tends to $0$ as $\Delta \rightarrow \infty$.
\end{theorem}

\renewcommand{\vol}{\mu}

Our proof of Theorem~\ref{th:main} goes as follows.  Take any small constant $\eps>0$ and assume that $m\geq m_0(\eps)$ is sufficiently large. Then take any number $R>0$. 
Take a lattice $\mathcal{L}\subseteq \CG$ of girth at least $16R$, which exists by Lemma~\ref{lm:cover}. Let $M:=\mathbb{H}^{m}/\mathcal{L}$ with measure $\vol=\mu_M$ and distance $\dist_M$.
 Consider $X\sim \mathrm{Po}_\lambda(M)$, that is, the Poisson process $X\subseteq M$ with intensity $\lambda \, \vol$ for some choice of $\lambda\in (0,\infty)$.
Once we obtain $X$, we will delete some bad points from $X$ to obtain $Y$, and consider $G_{Y}$, the graph with the vertex set $Y$ and the edge set 
\beq{eq:G}
 E(G_Y):=\{\,\{x,y\}: x,y\in Y \text{ and } 0<\dist_M(x,y)\le 2R\}.
 \eeq
Based on our choice of $\lambda$ and the removal process, the graph $G_{Y}$ will satisfy the condition in Theorem~\ref{thm: indep}. Thus Theorem~\ref{thm: indep} yields an independent set in $G_{Y}$, which corresponds to an $R$-packing in $M$. This gives a lower bound on $\D_R(M)$ which is also a lower bound on $\D_R(\I H^m)$ by Lemma~\ref{lm:Transfer}.

Let us give the details of the proof (for finding $Y$ after $\C L$ and $M$ have been defined). We will use asymptotic notation, like $O(1)$, with respect to $m\to\infty$ for constants that can be chosen independently of~$R$.

Given $m$ and $R$, consider the following function of $x\in (0,R]$:
\begin{equation}\label{eq:gamma}
 \gamma(x):=\left\{ \begin{array}{ll} m\cdot \tanh^2{(x/2)} -50\tanh^2(2R)\cdot \left(\ln(m) + \ln\ln\left(\frac{\sinh(2R)}{\sinh x}\right)\right),& \mbox{if $R <m$,}\\
 \cosh^2{(x/2)} - m\ln(\cosh(2R)),&\mbox{if $R\ge m$}.
\end{array}\right. 
\eeq
 Note that, by the monotonicity of $\sinh$, the double logarithm in~\eqref{eq:gamma} is well-defined.
 We define $\tau = \tau(m,R)$ to be the number between $0$ and $R$ satisfying $\gamma(\tau)=0$. 

Let us show that such $\tau$ uniquely exists (if $m$ is sufficiently large). Note that $\gamma(x)$ is a strictly increasing and continuous function of $x\in (0,R]$. So it is enough to show that it assumes both positive and negative values. If $R\ge m$ then we have rather generously that
$\gamma(R)\ge \me^{R}/4-2Rm>0$ and $\gamma(\ln m)\le m-m^2<0$, as desired. Moreover, we have $\tau > \ln m$ in this case.
So let us consider the case $R<m$.
 Using $0<\tanh(2R)<4\tanh (R/2)$, we obtain that 
 $$
 \gamma(R)\ge \tanh^2(R/2)\left(m -O(\ln m)\right)>0. 
 $$
 On the other hand, if we let $x\to 0$ then $\gamma(x)$ tends to $-\infty$. Thus $\tau\in (0,R]$ satisfying $\gamma(\tau)=0$ exists and is unique.

\begin{claim}\label{cl: basic}
If 
%$m\to\infty$ and 
$R<m$, then the following holds:
    $$\frac{\sinh (2R)}{\sinh \tau} = \Theta\left(\cosh(2R)\sqrt{\frac{m}{\ln(m)}}\right).$$
\end{claim}
\bpf Observe that $\tau=o(1)$ for otherwise $R\ge \tau\ge \Omega(1)$ and $\gamma(\tau)=\Omega(m)$ cannot be~0. 

Let $\kappa:=\ln m + \ln\ln\left(\frac{\sinh(2R)}{\sinh \tau}\right)$. By the monotonicity of $\sinh$ and by $\sinh 2x\ge 2\sinh x$ it holds that, for example,
\begin{equation}\label{eq:kappa}
\ln\left(\frac{\sinh(2R)}{\sinh x}\right) \ge\ln\left(\frac{\sinh(2R)}{\sinh R}\right) \ge \ln 2\ge \frac1{\sqrt m},\quad\mbox{for any $0<x\le R$}.
 %-\frac12\ln m\le \ln\ln 2\le \ln\ln\left(\frac{\sinh(2R)}{\sinh(R)}\right)\le \ln\ln\left(\frac{\sinh(2R)}{\sinh \tau}\right)\le \ln\ln(\sinh(2R))\le 2\ln m,
 %\quad\mbox{for any $0<x\le R$}.
 \end{equation}
 Thus $\kappa\ge \frac12\ln m$. 

It is enough to show that $\kappa=O(\ln m)$. Indeed, then
 $c:=\frac{\tanh^2(\tau)}{\tanh^2 (\tau/2)}$ satisfies $1\leq c\leq 4$ and we have
\begin{align*}
    \frac{\sinh^2(2R)}{\sinh^2(\tau) } 
&
%=  \frac{\cosh^2(2R)}{\cosh^2(\tau)}  \cdot \frac{\tanh^2(2R)}{\tanh^2(\tau)}
=  \frac{\cosh^2(2R)}{\cosh^2(\tau)}  \cdot \frac{\tanh^2(2R)}{c\tanh^2(\tau/2)}
= 
\frac{\cosh^2(2R)}{\cosh^2(\tau)} \cdot \frac{m}{50c\kappa} =  \frac{\cosh^2(2R)}{\cosh^2(\tau)} \cdot \frac{m}{\Theta(\ln m)}.
\end{align*}
 from which the claim follows by $\cosh \tau=1+o(1)=\Theta(1)$.

%as then $0=\gamma(\tau)=\Theta(m\tau^2)-\Theta(\tanh^2(2R)\ln m)$ from which the claim follows. 

The proof of $\kappa=O(\ln m)$ is routine except we have to be careful to rule out the possibility that $\tau$ is extremely small relative to $m$ and~$R$. Suppose on the contrary that $\kappa\not=O(\ln m)$. Hence $\ln(\frac{\sinh(2R)}{\tau}) = \omega(m)$. Then 
%necessarily $\tau=o(1)$ thus 
the identity $\gamma(\tau)=0$ implies by $\tau=o(1)$ that
 \beq{eq:Tau2}
 m\tau^2=\Theta\left(\frac{\sinh^2(2R)}{1+\sinh^2(2R)} \ln\ln\left(\frac{\sinh(
2R)}{\tau}\right)\right)
 \eeq
  It follows that $R=o(1)$ as otherwise $\tau\ge m^{-1/2+o(1)}$ by~\eqref{eq:Tau2}, contradicting $\ln(\frac{\sinh(2R)}{\tau}) = \omega(m)$. We obtain from~\eqref{eq:Tau2} that $m=\Theta\left((R/\tau)^2 \ln\ln(R/\tau)\right)$. Since $m\to\infty$, it follows that $R/\tau\to\infty$ and thus $R/\tau=\Theta\left(\sqrt{m/\ln\ln m}\right)$, which again contradicts our assumption $\kappa\neq O(\ln m)$ .\epf

Now we define the parameters $\Delta$ and $\lambda$ as follows:
$$\Delta := \frac{1}{m^4} \frac{\vol(\Ball_{2R})}{ \vol(\Ball_{\tau})}
\text{ and } \lambda := \frac{\Delta}{\vol(\Ball_{2R})} .$$
Recall that $\vol=\mu_M$ denotes the measure on $M$ and that,  for every $0\le r\le 4R$, it assigns the same measure to an $r$-ball in $M$ as the hyperbolic measure $\HMeasure$ assigns to an $r$-ball in $\I H^m$ by the girth assumption on $\CL$.

%In fact, as $\tau\leq R$, for either case, $\Delta$ is at least $2^{m+o(m)}$. We can further prove more accurate estimate on $\Delta$ as follows:
Let us show that \begin{align}\label{eq: lnDelta}
 \ln \Delta = (1+o(1)) m\ln (\sqrt{m} \cosh(2R)).
\end{align}
Indeed, if $R<m$ then \eqref{eq: lnDelta} is a consequence of Lemma~\ref{lem: vol ratio}  and Claim~\ref{cl: basic}, so suppose that $R\geq m$. Then we have $\tau> \ln m$ and $\sinh \tau \leq \cosh \tau \leq 2\cosh^2 (\tau/2)$, thus \eqref{eq:gamma} implies 
$$
0< \ln (\sinh \tau) = O(\ln R + \ln m) = o(\ln (\sinh (2R))).
$$
Thus
Lemma~\ref{lem: vol ratio} implies  that 
\begin{eqnarray*}
 \ln \Delta &=&  (m+O(1)) (\ln (\sinh (2R)) - \ln (\sinh \tau) ) +O(\ln m)\\
 &=&  (1+o(1))m \ln (\sinh (2R)), 
%&=& (1+o(1)) m\ln (\sqrt{m}\cosh (2R)),
\end{eqnarray*}
giving \eqref{eq: lnDelta} since $|\sinh (2R)- \cosh (2R)| \le 1$ and 
$ \ln (\sqrt{m})= o(\ln (\cosh (2R)))$.

In particular, we have that $\Delta\to\infty$ as $m\to\infty$.

With these choices, we have the following lemma. Recall that $G_{Y}$ is the graph on $Y$ whose edge set is defined by~\eqref{eq:G}.
%, where $G(Y,R)$ is the graph on the vertex set $Y$ and edge set $\{xy\in Y: y\in \Ball_{2R}(x)\}$.

\begin{lemma}\label{lem: discretization}
There exists $Y\subseteq M$ such that 
$$ |Y| \geq \left(1-\frac{1}{m}\right)\frac{\Delta}{\vol(\Ball_{2R})}\,\vol(M) $$
and the graph $G_{Y}$ satisfies that
$$\Delta(G_Y)\leq \Delta + \Delta^{2/3} \text{ and } \Delta_2(G_Y)\leq \Delta\, (\ln \Delta)^{-10}.$$
\end{lemma}

With this lemma and Theorem~\ref{thm: indep}, we obtain that there is an $R$-packing in $M$ of size at least 
$$
 (1-o(1))\frac{|Y|\ln \left(\Delta+ \Delta^{2/3}\right)}{\Delta+ \Delta^{2/3}}\ge (1-1/m-o(1))\frac{\vol(M) \ln \Delta}{\vol(\Ball_{2R})}.$$
Thus \eqref{eq: lnDelta} implies  that, as $m\to\infty$,
$$\D_R(\I H^m)\ge \D_R(M)\geq (1-o(1))\ln \Delta \cdot \frac{\vol(\Ball_{R})}{\vol(\Ball_{2R})}
\geq  (1-\eps )m \ln (\sqrt{m} \cosh(2R)) \frac{\vol(\Ball_{R})}{\vol(\Ball_{2R})}.$$
This proves~\eqref{eq:main} (that is, Theorem~\ref{th:main}). Thus it remains to prove Lemma~\ref{lem: discretization}.

\bpf[Proof of Lemma~\ref{lem: discretization}.]
We follow the argument in Section 2 of \cite{CamposJenssenMichelenSahasrabudhe}.
In order to prove Lemma~\ref{lem: discretization}, we sample a Poisson point process $X\subseteq M$ with the intensity $\lambda\,\vol$, 
and obtain a desired set $Y\subseteq X$ by removing points $x\in X$ which satisfy at least one of the following two conditions:
\begin{align}\label{eq: bad}
    |X\cap \Ball_{2R}(x)|\geq \Delta + \Delta^{2/3} \text{,\ \  or\ \  $\exists\, y\in X$\ \  } |X\cap \Ball_{2R}(x)\cap \Ball_{2R}(y)| \geq \Delta\, (\ln \Delta)^{-10}.
\end{align}

We will show that $Y$ has almost the same size as $X$. 
The following \emph{identity of Mecke} will be useful: for any bounded measurable set $\Lambda \subseteq M$ and events $(A_x)_{x\in \Lambda}$ we have
\begin{align}\label{eq: Mecke}
    \mathbb{E}\left| \{ x\in X\cap \Lambda: A_x \text{ holds for } X\}\right| = \lambda \int_{\Lambda} \mathbb{P}[\,A_x \text{ holds for } X\cup \{x\}\,] \dd\vol (x).
\end{align}
Also, the following tail bound for a Poisson random variable $Z$ will be used:
\begin{align}\label{eq: tail bound}
    \mathbb{P}[\,Z -\mathbb{E}Z \geq t\,\mathbb{E}Z\,] \leq \exp\left(- \min\{t,t^2\}\frac{\mathbb{E}Z}{3} \right).
\end{align}
First, we show the following claim stating that, on average, only a small fraction of $x\in X$ satisfies the first bad condition in \eqref{eq: bad}.

\begin{claim}\label{cl: first bad}
Let $X\sim \mathrm{Po}_{\lambda}(M)$. Then 
$$\mathbb{E}|\{x\in X: |X\cap \Ball_{2R}(x)|\geq \Delta + \Delta^{2/3} \}| \leq \frac{1}{m^2}\, \mathbb{E}|X|.$$
\end{claim}
\bpf Take any $x\in M$.
Recall that $|X\cap \Ball_{2R}(x)|$ is a Poisson random variable of mean $\lambda \vol(\Ball_{2R}) = \Delta$. Thus using \eqref{eq: tail bound}, we have 
$$ \mathbb{P}[\, |X\cap \Ball_{2R}(x)| \geq \Delta+ \Delta^{2/3}-1\,] \leq \exp\left(-\frac{1}{4}\Delta^{1/3}\right) \leq m^{-2}.$$
Thus we have by Mecke's identity that
$$\mathbb{E}\, \left|\{x\in X: |X\cap \Ball_{2R}(x)|\geq  \Delta+\Delta^{2/3}\}\right|=
\lambda \int_{M} \mathbb{P}[\,|X\cap \Ball_{2R}(x)|\geq  \Delta+\Delta^{2/3}-1\,]\dd\vol(x)\le \frac{\lambda\, \vol(M)}{m^2},
$$
 giving the desired by $\mathbb{E} |X|=\lambda\, \vol(M)$.
\epf

Before considering the second bad condition in~\eqref{eq:  bad}, we prove the following claim.

\begin{claim}\label{cl: covolume}
For any $x,y\in M$ at distance $\dist_M(x,y)\ge \tau$, it holds that
\[\lambda \,\vol( \Ball_{2R}(x)\cap \Ball_{2R}(y))\leq \Delta\, (\ln \Delta)^{-15}.
\]
\end{claim}
\bpf
Note that we only have to consider the case where $\tau$ is at most $4R$, as otherwise $\Ball_{2R}(x)\cap \Ball_{2R}(y)$ is empty.
As the lattice $\CL$ has girth at least $16R$, Lemma~\ref{lm:Intersection} with $r=2R$ applies to the points $x,y$ in $M$ (instead of $\I H^m$). Hence we obtain that the following where $\sigma =\sigma(2R,\tau)\in (0,2R)$ was defined in~\eqref{eq:Intersection}:
\begin{align}\label{eq: expression}
\begin{split}
     \lambda\, \vol( \Ball_{2R}(x)\cap \Ball_{2R}(y))
  &\leq  \lambda\,  \vol(\Ball_{\sigma}) \leq  \Delta\, \frac{\vol(\Ball_{\sigma})}{\vol(\Ball_{2R})} \\
  &\leq 
   \Delta\, \frac{ \sinh^{m-1}( \sigma)}{\sinh^{m-1}(2R)} =  \Delta\,\left(1 - \coth^2(2R) \tanh^2\left(\tau/2\right)\right)^{(m-1)/2}.
\end{split}
\end{align}
Here, the penultimate inequality is from Lemma~\ref{lem: vol ratio}.

Recall that $\tau$ satisfies $\gamma(\tau)=0$ where $\gamma$ was defined in \eqref{eq:gamma}. If $R<m$, then the final expression in~\eqref{eq: expression} becomes
\begin{align*}
 \Delta \left(1 - 50\cdot \frac{\ln m + \ln\ln(\frac{\sinh 2R}{\sinh \tau} )}{m}\right)^{(m-1)/2}
   \leq  \Delta\, \frac{1}{m^{20}}\cdot \ln\left(\frac{\sinh 2R}{\sinh \tau} \right)^{-20}  \leq  \Delta\, (\ln \Delta)^{-15},
\end{align*}
as we want. Here, the final inequality holds as Lemma~\ref{lem: vol ratio} implies $ \ln \Delta \leq  m \ln(\frac{\sinh 2R}{\sinh \tau})$.
If $R>m$, then 
as $\coth^2(2R)\tanh^2(\tau/2) \ge \tanh^2(\tau/2)
%=  \left(1+\frac{1}{\cosh^2(2R)-1}\right)\left(1- \frac{1}{\cosh^2 (\tau/2)}\right) 
= 1 - \frac{1}{\cosh^2(\tau/2)}$, 
the final term in \eqref{eq: expression} is bounded from above by
\begin{align*}
 \Delta \left( \frac{1}{\cosh^2{(\tau/2)}}\right)^{(m-1)/2} 
    &= \Delta \left(\frac{1}{m \ln (\cosh(2R))} \right)^{(m-1)/2} \\
    & \leq \Delta\,  (\ln \Delta)^{-(m-1)/4} \leq \Delta\,  (\ln \Delta)^{-15},
\end{align*}
where the penultimate inequality follows from \eqref{eq: lnDelta}.
This proves the claim.
\epf

Now we bound the number of points satisfying the second bad condition in \eqref{eq: bad}. 
\begin{claim}\label{cl: second bad}
Let $X\sim \mathrm{Po}_{\lambda}(M)$ and put $\eta := (\ln \Delta)^{-10}$. Then we have that
$$\mathbb{E}|\{x\in X: |X\cap \Ball_{2R}(x)\cap \Ball_{2R}(y)|\geq \eta \Delta \text{ for some } y\in X\}| \leq \frac{1}{2m}\, \mathbb{E}|X|.$$
\end{claim}
\bpf
Take any $x\in M$. For $y\in M$, let $I_{x,y}:= |X\cap \Ball_{2R}(x)\cap \Ball_{2R}(y)|$.
%Note that $I_{x,y}$ is a Poisson random variable.
Using Markov's inequality, we have 
$$\mathbb{P}[\,\exists y\in X: I_{x,y} \geq \eta \Delta -1\,] 
\leq \mathbb{E}\,|\Ball_{\tau}(x)\cap X| + \mathbb{E}\,|\{y\in X\setminus \Ball_{\tau}(x) \mid I_{x,y} \geq \eta \Delta -1\} |.$$
We bound each of these two terms.

For the first term, using Lemma~\ref{lem: vol ratio} and the definition of $\Delta$, we have
$$\mathbb{E}\,|\Ball_{\tau}(x)\cap X|\le \lambda\, \vol(\Ball_{\tau}(x)) \leq 
\frac{\Delta\,\vol(\Ball_{\tau})}{ \vol(\Ball_{2R})}
= \frac{1}{m^4}. $$

For the second term, we only need to consider $y\in \Ball_{4R}(x)$ as otherwise $I_{x,y}=0<\eta \Delta-1$. (Note that $\eta\Delta\to\infty$.)
Using \eqref{eq: Mecke} and Markov's inequality, the second term is 
\begin{align*} 
\lambda\, &\int_{M\setminus\Ball_{\tau}(x)} \mathbb{P}[\,I_{x,y} \geq \eta \Delta -2 \,] \dd\vol(y) 
  = \lambda\, \int_{\Ball_{4R}(x)\setminus\Ball_{\tau}(x)} \mathbb{P}[\,I_{x,y} \geq \eta \Delta -2\, ] \dd\vol(y) \\
&\leq \lambda\,  \vol(B_{4R})\ 
%\cdot\ \mathbb{E}|X\cap \Ball_{4R}(x)|\ \cdot\ 
\sup_{y\in \Ball_{4R}(x) \setminus \Ball_{\tau}(x)} \mathbb{P}[\,I_{x,y} \geq \eta \Delta -2\,].
\end{align*}
Using Claim~\ref{cl: covolume}, we have for $y$ with $\dist_M(y,x)\ge \tau$ that
$$\mathbb{E} I_{x,y}= \lambda\, \vol( \Ball_{2R}(x)\cap \Ball_{2R}(y))\leq \Delta\, (\ln \Delta)^{-15}.$$
As $I_{x,y}$ is a Poisson random variable, we can apply \eqref{eq: tail bound} to obtain that, rather roughly,
\begin{align*}
    \lambda\,  \vol(B_{4R})\
    %\mathbb{E}|X\cap \Ball_{4R}(x)|\ 
    \mathbb{P}[\,I_{x,y}\geq \eta \Delta -2\,] &\leq 
\Delta\, \frac{\vol(\Ball_{4R})}{\vol(\Ball_{2R})} \exp\left(-\Delta (\ln\Delta)^{-15}\right) < \frac{1}{m^2}.
%\\ &\leq \frac{\Delta}{m^4} \frac{\vol(\Ball_{4R})}{\vol(\Ball_{2R})}  \exp(-\Delta (\ln\Delta)^{-15}) \leq \frac{1}{m^2}.
\end{align*}
The last inequality follow by observing by Lemma~\ref{lem: vol ratio} and \eqref{eq: lnDelta} that
$$
 \frac{\vol(\Ball_{4R})}{\vol(\Ball_{2R})}\le \left(\frac{\sinh(4R)}{\sinh(2R)}\right)^m\le \me^{2Rm}  \leq \frac{\exp\left(\Delta (\ln\Delta)^{-15}\right)}{\Delta\, m^2}.
 $$

%The final inequality holds as Lemma~\ref{lem: vol ratio} implies 
%$\frac{\vol(\Ball_{4R})}{\vol(\Ball_{2R})} \leq 
%\left(\frac{\sinh (4R)}{\sinh(2R)}\right)^{m-1} \leq (2\cosh 2R)^{m-1} \leq %\Delta^{2}.$

Thus we have shown that, for every $x\in X$, 
$$
\mathbb{P}[\,\exists y\in X: I_{x,y} \geq \eta \Delta -1\,] \le \frac1{m^4} + \frac1{m^2}<\frac1{2m}.
$$
By Mecke's idenity, we conclude that the expected number of $x\in X$ satisfying the second bad condition in~\eqref{eq: bad} is at most $\frac 1{2m}\lambda\,\vol(M) $, proving the claim.
 \epf

Now we prove Lemma~\ref{lem: discretization}.
Let $X$ be obtained from the Poisson point process in $M$ with intensity~$\lambda\, \vol$.
Let $S_1, S_2\subseteq X$ be the points that satisfy the first and second bad properties in \eqref{eq: bad}. Let $Y:=X\setminus (S_1\cup S_2)$.
Then the previous claims imply that 
$$\mathbb{E}|Y| \geq \mathbb{E}|X|- \mathbb{E}|S_1|-\mathbb{E}|S_2| \geq \left(1-\frac1{m^2}-\frac{1}{2m}\right)\mathbb{E}|X| \geq 
\left(1-\frac{1}{m}\right)\frac{\Delta}{\vol(\Ball_{2R})}\, \vol(M).$$
There is an outcome $X$ such that $|Y|$ is at least its expected value, finishing the proof of Lemma~\ref{lem: discretization}.\epf

\section*{Acknowledgements}

The authors would like to thank Lewis Bowen for the helpful discussions and the anonymous referees for their useful comments.

%\newpage

%\bibliography{bibexport}
%	\bibliography{oleg,sets,misc,ramsey,enum,number,posets,sat,ex,matroid,design,random,graph,general,geometry,algorithm,Analysis,limits}


\hide{
\section*{Macros}

	$\D_R(V)$ \verb$\D_R(V)$: $R$-packing density in metric space $V$
	
	$\I H^m$ \verb$\I H^m$: hyperbolic space

	$\CG$ \verb$\CG$: group of isometries of $\I H^m$
	
	$\CL$ \verb$\CL$: lattice
	
	$\CO$ \verb$\CO$: fixed point in $\I H^m$
	
	$\dist_m$ \verb$\dist_m$: distance on $\I H^m$
	
    $\HMeasure$ \verb$\HMeasure$: measure on $\I H^m$

	$\X,\Y$ \verb$\X,\Y$: sphere packings in $\I H^m$
	
	$\dist_R$ \verb$\dist_R$: distance on $S_R$
	
$\dist_M,\mu_M$ \verb$\dist_M,\mu_M$: distance and measure on $M$

	$\D_R(\I H^m)$ \verb$\D_R(\I H^m)$: Bowen--Radin  radius-$R$ packing density
	
	$\Orbit(\X)$ \verb$\Orbit(\X)$: orbit of a packing $X$ in $S_R$
	
	$\Ball_R(\X)$ \verb$\Ball_R(\X)$: $R$-ball around $\X$
}	
	
\end{document}